\documentclass[10pt]{article}

\newtheorem{theorem}{Theorem}
\newtheorem{lemma}{Lemma}

\date{}

\title{The Limiting Distribution of the Trace of a Random Plane Partition}
\bigskip

\author{Emil P. Kamenov \footnote{The author was partially supported by the
*MRTN-CT-2004-511953* project grant carried out by the Alfr\'{e}d R\'{e}nyi
Institute of Mathematics in the framework of the European Community's Program:
"Structuring the European Research Area".}\\
    \small\it Sofia University\\
    \small and\\
    Ljuben R. Mutafchiev\footnote{Partial support given by the National Science
    Fund of the Bulgarian Ministry of Education and Science, grant No. VU-MI-105/2005.}\\
    \small\it American University in Bulgaria and\\
    \small\it Institute of Mathematics and Informatics of the Bulgarian Academy of Sciences}

\begin{document}
\maketitle

\begin{abstract}
We study the asymptotic behaviour of the trace (the sum of the
diagonal parts) $\tau_n = \tau_n( \omega)$ of a plane partition
$\omega$ of the positive integer $n$, assuming that $\omega$ is
chosen uniformly at random from the set of all such partitions. We
prove that $(\tau_n -c_0 n^{2/3})/c_1 n^{1/3} \log^{1/2}n$
converges weakly, as $n\rightarrow \infty$, to the standard normal
distribution, where $c_0 = \zeta(2) / [2\zeta(3)]^{2/3}$, $c_1 =
\sqrt{1/3} / [2\zeta(3)] ^{1/3}$ and $\zeta(s) = \sum_{j=1}
^\infty j^{-s}$.
\end{abstract}

\def\thefootnote{}

\section{Introduction}

Properties of various kinds of partitions are often studied using
bivariate generating functions of the following type:

\begin{equation}\label{bivariate}
G(u,x;\{a_j\}_{j\geq 1}) = \prod_{j=1}^\infty(1-ux^j)^{-a_j} = 1 +
\sum_{n, m \geq 1}  Q(m,n;\{a_j\}_{j\geq 1} )\ u^m \ x^n.
\end{equation}
Here $u$ is finite, $|x|<1$ and $\{a_j\}_{j\geq 1}$ is a given
sequence of non-negative numbers. A combinatorial interpretation
of (\ref{bivariate}) for integer-valued sequences $\{a_j\}_{j\geq
1}$, is obtained expanding the geometric progressions in the
left-hand side. Setting
\begin{displaymath}
\Lambda = \bigcup_{k=1}^\infty \bigcup_{j=1}^{a_k} \{k\} =
\{\lambda_1,\lambda_2,\ldots \},
\end{displaymath}
one can show that the coefficient $Q(m,n;\{a_j\}_{j\geq 1})$
equals the number of solutions in non-negative integers $z_j$ of
the equations
\begin{displaymath}
\left\{ \matrix{\sum\limits_{\lambda_j\in\Lambda} \lambda_j z_j =
n,  \cr \cr \sum\limits_j z_j=m.}\right.
\end{displaymath}

Furthermore, a specific sequence of integers $\{a_j\}_{j\geq 1}$
may define a particular class of partitions, so that the $n$th
coefficient $Q(n;\{a_j\}_{j\geq 1})$ in the power series expansion
of
\begin{equation}\label{univariate}
G(1,x;\{a_j\}_{j\geq 1}) = \prod_{j=1}^\infty(1- x^j)^{-a_j} = 1 +
\sum_{n  \geq 1}  Q(n;\{a_j\}_{j\geq 1} )  x^n
\end{equation}
counts the number of partitions of $n$ belonging to that class
(see e.g. \cite[Chap. 6]{Andrews}). A fundamental problem in the
theory of partitions is to determine the asymptotic behaviour of
$Q(n;\{a_j\}_{j\geq 1})$ as $n\rightarrow \infty$. In its most
general setting this problem was studied by Meinardus
\cite{Meinardus} who introduced a set of assumptions on the
sequence $\{a_j\}_{j\geq 1}$ and obtained an expression for the
leading term in the asymptotic expansion of this coefficient. Two
important classes of partitions are covered by Meinardus' formula:
(i) $a_j =1$ and (ii) $a_j = j, j = 1, 2, \ldots$. In the first
case $Q(n;\{a_j = 1 \}_{j\geq 1})$ equals the number of (linear)
integer partitions of $n$ (see \cite[Section 1.2]{Andrews}) and
Meinardus' asymptotic formula implies the famous result of Hardy
and Ramanujan \cite{Hardy} for the number of such partitions. In
the second case $Q(n;\{a_j = j \}_{j\geq 1})$ equals the number of
plane partitions of $n$, whose asymptotic expression was
previously obtained by Wright \cite{Wright}.

A basic restriction in Meinardus' scheme states that the
Dirichlet's series
\begin{equation}\label{deots}
D(s) = \sum_{j =1}^\infty a_j j^{-s},
\end{equation}
generated by the sequence $\{a_j\}_{j\geq 1}$ has to converge in
the half-plane $\Re e\: s > \alpha >0$. It turns out that the
asymptotic behaviour of the coefficients $Q(m,n;\{a_j\}_{j\geq
1})$ in the bivariate power series expansion (\ref{bivariate})
strongly depends on the value of the parameter $\alpha$.
Haselgrove and Temperley \cite{Haselgrove} studied the case when
$\alpha < 2$, $n \rightarrow \infty$ and $m$ becomes large with
$n$ at a specific rate. They applied the classical method due to
Hardy and Ramanujan \cite{Hardy} and obtained a result of the form
of a local limit theorem for the ratio $Q(m,n;\{a_j\}_{j\geq 1} )
/ Q(n;\{a_j\}_{j\geq 1})$. Their result established convergence to
a non-Gaussian distribution. Haselgrove and Temperley
\cite[Section 3]{Haselgrove} also conjectured that the Gaussian
law would appear if the Dirichlet's series parameter $\alpha$ is
not less than 2, however, a formal proof is still lacking.

In this paper we consider the sequence $a_j  = j, j = 1, 2,
\ldots$, and focused on plane partitions of positive integer $n$.
A plane partition $\omega$ of $n$ is a representation
\begin{displaymath}
n = \sum_{i,j \geq 1} \omega_{i,j},
\end{displaymath}
in which the array $\omega  = (\omega_{i,j})_{i ,j \geq 1}$ of
non-negative integer entries is such that $\omega_{i,j} \geq
\omega_{i+1,j}$ and $\omega_{i,j} \geq \omega_{i,j+1}$. We may
also assume that $\omega$ occupies the first quadrant of the
coordinate system $iOj$. As an illustration the following is the
plane partition of $n=8$:
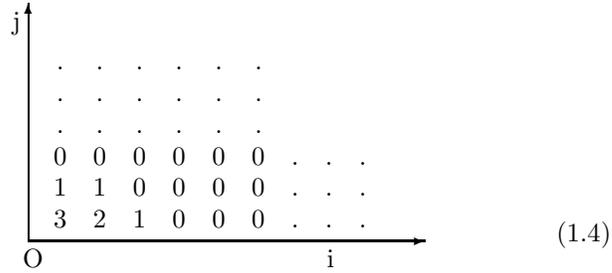
\begin{equation}\label{eq.example}
\begin{picture}(100,100)
\put(10,35){$\begin{matrix}{.&.&.&.&.&.&\cr.&.&.&.&.&.&\cr.&.&.&.&.&.&\cr
0&0&0&0&0&0&.&.&.\cr 1&1&0&0&0&0&.&.&.
\cr3&2&1&0&0&0&.&.&.}\end{matrix}$ } \put(2,0){\vector(0,1){90}}
\put(2,0){\vector(1,0){150}} \put(115,-10){i} \put(-4,80){j}
\put(0,-10){O}
\end{picture}
\end{equation}
For the sake of brevity the zeroes in (\ref{eq.example}) are
deleted, so that the abbreviation $$ \matrix{1&1 \cr 3&2&1\cr} $$
presents the plane partition of $8$ in a shorter way. It seems
that MacMahon \cite{Macmahon} was the first who suggested the idea
of a plane partition. Various properties of plane partitions and
their applications to combinatorics, algebra and analysis of
algorithms may be found in \cite[Chap. 11]{Andrews}, \cite[Chap.
12]{Nijenhuis} and \cite[Chap. 7]{Stanley3}.

Further, for the sake of brevity, we also let

\begin{equation}\label{brevity}
Q(n) = Q(n;\{a_j = j \}_{j\geq 1}), \hspace{.8cm} G(u ,x) =
\prod_{j=1}^\infty(1-  u x^j)^{-j}.
\end{equation}
In terms of notations (\ref{brevity}) the generating function
(\ref{univariate}) becomes
\begin{equation}\label{brevg}
G(1 ,x) = \prod_{j=1}^\infty(1- x^j)^{-j} = 1 + \sum_{n  \geq 1}
Q(n)  x^n.
\end{equation}
It is also easily seen that the sequence $a_j  = j, j = 1, 2,
\ldots$, implies that
\begin{equation}\label{zeta}
D(s) = \zeta(s-1) =  \sum_{j =1}^\infty j^{-(s-1)}
\end{equation}
(see (\ref{deots})). Therefore, in this case we have $\alpha = 2$.

The main goal of our paper is to prove Haselgrove and Temperley's
conjecture \cite[Section 3]{Haselgrove} that the Gaussian law
determines the asymptotic behaviour of $Q(m ,n) = Q(m, n;\{a_j = j
\}_{j\geq 1})$ if $m$ becomes large with $n$ at a specific rate.

It turns out that this problem has an interesting probabilistic
interpretation. If we introduce the uniform probability measure $P
= P_n$ on the set of all plane partitions of $n$ assuming that the
probability $1/Q(n)$ is assigned to each plane partition $\omega$,
then each conceivable numerical characteristic of $\omega$ becomes
a random variable. An exact combinatorial expression for the
numbers $Q(n)$ does not exist. As it was previously mentioned
their asymptotic was determined by Wright \cite{Wright} and
subsequently his result was confirmed by Meinardus' general
theorem \cite{Meinardus}. It was shown that
\begin{equation}\label{eq.wright}
Q(n)\sim \frac{[\zeta(3)]^{7/36}}{2^{11/36}3^{1/2}\pi^{1/2}}\
n^{-25/36} \exp \left\{3[\zeta(3)]^{1/3} (n/2)^{2/3} + 2c
\right\},
\end{equation}
where
\begin{displaymath}
  c = \int_0^\infty \frac{y \log y}{e^{2 \pi y}-1}d y.
\end{displaymath}
(In the statement of Wright's main result \cite[p.179]{Wright} the
constant $3^{1/2}$ in the denominator of (\ref{eq.wright}) is
missing, however, it is included in his final result at the end of
the proof of his theorem on p.189. Further results on the
asymptotic expansion of $Q(n)$ can be also found in
\cite{Almkvist, AlmkvistII}.)

We will consider here the trace $\tau_n$ of a partition $\omega$
defined as the sum of its diagonal parts:
\begin{displaymath}
\tau_n = \sum_{j=1}^n \omega_{j,j}.
\end{displaymath}
To study the asymptotic behaviour of $\tau_n$ as $n\rightarrow
\infty$ we use the following generating function identity (see
\cite{Stanley} or \cite[Chap. 11, Problem 5]{Andrews}) :
\begin{eqnarray}\label{eq.stanley}
\nonumber G(u,x) = 1 + \sum_{n=1}^\infty Q(n)\ x^n \sum_{m=1}^n
P(\tau_n=m)\ u^m = \\ = 1 + \sum_{n=1}^\infty Q(n)\ \varphi_n(u)\
x^n = \prod_{j=1}^\infty(1-ux^j)^{-j},
\end{eqnarray}
where
\begin{equation}\label{phi}
\varphi_n(u) = \sum_{m=1}^n P(\tau_n=m)u^m
\end{equation}
is the probability generating function of the trace $\tau_n$ and
$G(u,x)$ is the generating function defined by (\ref{brevity}). It
is now clear that if Haselgrove and Temperley's conjecture, that
we stated above, is valid, then the trace $\tau_n$ of a random
plane partition has to be asymptotically normal as $n\rightarrow
\infty$. The main result of this paper is the following limit
theorem.
\begin{theorem} For any real and finite $z$, we have
\begin{displaymath}
\lim_{n \rightarrow \infty}P\left( \frac{\tau_n - c_0n^{2/3}}{c_1
n^{1/3} \log^{1/2} n}\leq z \right) = \frac{1}{\sqrt{2\pi}}
\int\limits_{-\infty}^z e^{-y^2/2} dy,
\end{displaymath}
where
\begin{displaymath}
c_0 = \zeta(2) / \left[2\zeta(3)\right]^{2/3} = .916597104,$$
$$c_1 = \sqrt{1/3} / [2\zeta(3)] ^{1/3} = .430977269.
\end{displaymath}
\end{theorem}

This result confirms Haselgrove and Temperley's conjecture when
$\alpha = 2$ and $D(s) = \zeta(s-1)$ ( see (\ref{zeta})). We
believe that our method of proof can be utilized to study the
general case $\alpha \geq 2$.

Our study is also partially motivated by similar results obtained
for linear partitions of $n$. Erd\H{o}s and Lehner \cite{Erdos}
were apparently the first who have studied random linear
partitions using a probabilistic approach. As a matter of fact,
they showed that the number of summands, after an appropriate
normalization, converges weekly, as $n \to \infty$, to a random
variable having the extreme value distribution. Subsequent work by
a number of authors provides considerable information about the
structure of a "typical" linear partition of a large integer. (We
refer the reader e.g. to \cite{Szalay}, \cite{Szalay2},
\cite{Szalay3}, \cite{Erdos2}, \cite{Fristedt}, \cite{Pittel1},
\cite{Hwang}, \cite{Mutafchiev2} and the references therein.) To
this effect, Theorem 1 continues the study of partitions initiated
by the probabilistic approach.

The proof of our main result is based on a method developed by
Hayman \cite{Hayman}. To get an estimate  for the Cauchy integral
stemming from (\ref{eq.stanley}) we use the fact that the
generating function $G(1,x)$ defined by (\ref{brevg}) satisfies
Hayman's admissibility conditions in a neighborhood of its main
singularity $x=1$ and outside it (see Lemmas 1 and 2,
respectively). A relevant approach to problems related to other
characteristics of random plane partitions may be found in
\cite{Pittel2} and \cite{Mutafchiev3}.

We organize our paper as follows. Section 2 contains auxiliary
facts on the admissibility of the generating function $G(1,x)$
defined by (\ref{brevg}) and on the asymptotic behaviour of the
coefficients $Q(n)$ in its power series expansion. In Section 3 we
present the proof of Theorem 1.


\section{Preliminary Asymptotics}
\setcounter{equation}{0}

In order to get an asymptotical estimate for $G(1,x)$ around its
main singularity $x=1$, we first introduce the analytic scheme of
assumptions on the sequence $\{a_j\}_{j\geq 1}$ of non-negative
numbers, which is due to Meinardus \cite{Meinardus}. The following
three conditions must be satisfied:

(M1) The Dirichlet series $D(s)$ (see (\ref{deots})) converges in
the half-plane  $\Re e\ s>\alpha>0$, and can be analytically
continued into the half-plane  $\Re e\ s\geq -\alpha_0$, $\alpha_0
\in (0,1)$. In $\Re e\ s\geq -\alpha_0$, $D(s)$ is analytic except
for a simple pole at $s=\alpha$ with residue $A$.

(M2) There exists an absolute constant $\alpha_1 > 0$ such that
$D(s)= O(|\Im m\ s|^{\alpha_1})$ uniformly for $\Re e\ s \geq -
\alpha_0$ as $|\Im m \ s|\rightarrow \infty$.

(M3) Define $g(v)=\sum\limits_{j=1}^\infty a_j e^{-jv}$, where
$v=y + 2\pi i w$ and $y$ and $w$ are real numbers. If $|\arg
v|>\pi/4$ and $|w|\leq 1/2$, then $\Re e\ g(v)-g(y) \leq -\alpha_2
y^{-\chi}$ for sufficiently small $y$, where $\chi >0$ is an
arbitrary number, and $\alpha_2 >0$ is suitably chosen and may
depend on $\chi$.

We shall be concerned, in the first instance, with the behavior of
$G(1,x;\{a_j\}_{j\geq 1})$ (see \ref{univariate}) as $x$ becomes
close to 1.  Meinardus \cite{Meinardus} (see also Andrews
\cite[Lemma 6.1] {Andrews} proved that under the assumptions (M1)
and (M2)
\begin{equation}\label{eq.meinardus}
G(1,e^{-v};\{a_j\}_{j\geq 1}) = exp \left\{ A \Gamma(\alpha)\zeta
(\alpha +1)v^{-\alpha} - D(0) \log v + D'(0)+
O(y^{\alpha_0})\right\}
\end{equation}
as $y\rightarrow 0$ uniformly for $|\arg v|\leq \pi/4$ and $|w|
\leq 1/2$. (Here $\Gamma(\alpha)$ denotes Euler's gamma function
and $\log{(\cdot)}$ presents the main branch of the logarithmic
function satisfying $ \log {v}<0$ for $0<v<1$.)

Let us take now a sequence $\{r_n\}$ which, as
$n\rightarrow\infty$, satisfies
\begin{equation}\label{eq.defrn}
r_n = 1- \frac{\left[2\ \zeta(3)\right]^{1/3}}{n^{1/3}} +
\frac{\left[2\ \zeta(3)\right]^{2/3}} {2 n^{2/3}} -
\frac{\zeta(3)}{3\,n} + O(n^{-4/3}).
\end{equation}
For the sake of brevity we also set
\begin{equation}\label{eq.defb}
b(r)= \frac{6\ \zeta(3)}{(1-r)^4},
\end{equation}
where $0<r<1$. It is easy to check that (\ref{eq.defrn}) and
(\ref{eq.defb}) imply
\begin{equation}\label{eq.brn}
b(r_n)=\frac{3n^{4/3}}{[2\ \zeta (3)]^{1/3}} + O(n)
\end{equation}
as $n\rightarrow\infty$.

The next lemma suggests a tool that we shall subsequently use in
Section 3 to obtain the main term in our asymptotics.


\begin{lemma}
If $r_n$ satisfies (\ref{eq.defrn}) for large $n$, then
\begin{displaymath} G(1,r_ne^{i\theta})e^{-i\theta n } = G(1,r_n)e^{-\theta^2 b(r_n)/2}\left[1 +
O(1/\log^3 n)\right]
\end{displaymath}
as $n\rightarrow\infty$ uniformly for $|\theta|\leq \delta_n $,
where
\begin{equation}\label{eq.defdelta}
\delta_n = \frac{n^{-5/9}}{\log n}
\end{equation}
and $b(r_n)$ is determined by (\ref{eq.defb}).
\end{lemma}

\begin{proof} Our starting point here will be Meinardus' general asymptotic formula (\ref{eq.meinardus}).
We apply it for the sequence $a_j = j, j = 1,2, \ldots$. It is not
difficult to show that $A=1$ and $D(0)=-1/12$. Classical result on
the $\zeta$ function implies condition (M2) (see e.g.
\cite[Section 13.51 ]{Whittaker}). Therefore, for $v=y + 2\pi i
w$, we get
\begin{equation}\label{eq.c1}
G(1,e^{-v}) = exp\left\{\zeta(3) v^{-2} + \frac{1}{12}\ \log v +
D'(0) + O(y^{\alpha_0}) \right\}
\end{equation}
as $y\rightarrow 0$ uniformly for $|w|\leq 1/2$ and $|\arg v| \leq
\pi /4$. Setting
\begin{equation}
e^{-v} = r_n \ e^{i \theta},
\end{equation}
we see that $y=y_n=-\log{ r_n}$ and $w = -\theta / 2\pi$. The
asymptotic behaviour of $ -\log {r_n}$ can be determined with aid
of (\ref{eq.defrn}) as follows:
\begin{displaymath}
y_n= -\log{r_n} = \frac{\left[2\ \zeta(3)\right]^{1/3}}{n^{1/3}} -
\frac{\left[2\ \zeta(3)\right ] ^{2/3}}{2\, n^{2/3}} +
\frac{\zeta(3)}{3\,n}
\end{displaymath}
\begin{displaymath}
+{\left[2\ \zeta(3)\right ] ^{2/3} \over 2 n^{2/3}} - {2^{4/3}\
\left[2\ \zeta(3)\right ] ^{2/3} \left[\zeta(3)\right ] ^{1/3}
\over 4n} +\frac{2\, \zeta(3)}{3n} + O(n^{-4/3})
\end{displaymath}
\begin{equation}\label{eq.defyn}
={\left[2\ \zeta(3)\right]^{1/3}\over n^{1/3}} + O(n^{-4/3}),\
n\rightarrow \infty.
\end{equation}
Combining (\ref{eq.c1}) - (\ref{eq.defyn}), we observe that
\begin{equation}\label{eq.c2}
\frac{G(1,r_ne^{i\theta})}{G(1, r_n)}e^{-i\theta n}
\end{equation}
\begin{displaymath}
=\left({y_n - i\theta \over y_n}\right) ^{1/12} \exp
\left\{\zeta(3) [ (y_n - i\theta)^{-2} - y_n^{-2}] - i\theta n +
O\left( y_n^{\alpha_0}\right) \right\}.
\end{displaymath}
A Taylor's formula expansion for $|\theta| \leq \delta_n$ yields
\begin{displaymath}
(y_n  - i\theta)^{-2} - y_n^{-2} = 2i \theta y_n ^{-3} - 6
{\theta^2 \over 2} y_n^{-4} + O(|\theta|^3 y_n^{-5})
\end{displaymath}
\begin{equation}\label{eq.c3}
= 2i \theta y_n ^{-3}- 3 \theta^2  y_n^{-4} + O(\delta_n^3
y_n^{-5}).
\end{equation}
Using (\ref{eq.defyn}), we also get the following estimate for the
factor outside the exponent in (\ref{eq.c2}):
\begin{displaymath}
\left({y_n  - i\theta \over y_n}\right) ^{1/12} = \left\{ [2
\zeta(3)]^{1/3} - i \theta n^{1/3} + O(n^{-1}) \over  [2
\zeta(3)]^{1/3} + O(n^{-1}) \right\}^{1/12}
\end{displaymath}
\begin{equation}\label{eq.c4}
=1 + O(\delta_n n^{1/3}).
\end{equation}
Finally, we notice that (\ref{eq.defyn}) implies the bound
\begin{equation}\label{eq.c5}
y_n^{\alpha_0} = O(n^{-\alpha_0 /3}).
\end{equation}
Hence, inserting (\ref{eq.brn}),
(\ref{eq.defyn}),(\ref{eq.c3})-(\ref{eq.c5}) into (\ref{eq.c2}),
we obtain
\begin{displaymath}
{G(1,r_ne^{i\theta}) \over G(1, r_n)}e^{-i\theta n} = \left[ 1 +
O(\delta_n n^{1/3})\right]\exp\left\{\zeta(3) \left[{2i \theta
n\over 2 \zeta(3)[1+O(n^{-1})]^3}\right. \right.
\end{displaymath}
\begin{displaymath}
\left. \left. -{\theta^2 \over 2}{ 6n^{4/3} \over 2^{4/3}
[\zeta(3)]^{4/3}  [1+O(n^{-1})]^4} + O\left(\delta_n^3
y_n^{-5}\right)\right] -i \theta n + O\left(n^{-\alpha_0
/3}\right)\right\}
\end{displaymath}
\begin{displaymath}
= \left[ 1 + O(n^{-2/9}/\log{n})\right]\exp\left\{\zeta(3)
\left[{i \theta n\over \zeta(3)} \left(1+O(n^{-1})\right)
\right.\right.
\end{displaymath}
\begin{displaymath}
\left. \left. -\frac {\theta^2 }{ 2 \, \zeta(3)}
\left[b(r_n)+O(n)\right] \left(1+O(n^{-1})\right) + O(\delta_n^3
n^{5/3})\right] - i \theta n + O\left(n^{-\alpha_0 /3}\right)
\right\}
\end{displaymath}
\begin{displaymath}
= \left[ 1 + O(n^{-2/9}/\log{n})\right]\exp\left\{ i \theta n +
O(\delta_n) - \theta^2 b(r_n)/2 + O\left(n
\delta_n^2\right)\right.
\end{displaymath}
\begin{displaymath}
\left. +  O(\delta_n^2 b(r_n) n^{-1}) + O(\delta_n^3 n^{5/3}) - i
\theta n + O\left(n^{-\alpha_0 /3} \right) \right\}
\end{displaymath}
\begin{displaymath}
=\left[ 1 + O(n^{-2/9}/\log{n})\right]\exp\left\{ - \theta^2
b(r_n)/2 + O(1/\log^3{n}) \right\}
\end{displaymath}
\begin{displaymath}
= e^{- \theta^2 b(r_n)/2}\left[ 1 + O(1/\log^3{n})\right].
\end{displaymath}
This completes the proof.
\end{proof}


We also need another lemma that will establish a uniform estimate
for $G(e^{i T}, x)$, $|x|=r_n$ outside the range $-\delta_n <
\arg{x} < \delta_n$, if $T$ is real and suitably bounded.

\begin{lemma}
If $r_n$ and $\delta_n$ satisfy (\ref{eq.defrn}) and
(\ref{eq.defdelta}), respectively, and the function $T=T(n)$ is
such that
\begin{equation}\label{def.T}
T = T(n) = \Theta ( n^{-1/3} / \sqrt{ \log{n} })
\end{equation}
as $n\rightarrow \infty$, then there exist two positive constants
$\epsilon$ and $n_0$ so that
\begin{displaymath}
|G(e^{iT},r_ne^{i\theta})| \leq G(1,r_n) \exp
\left\{\varepsilon-2n^{2/9}/[2\zeta(3)]^{4/3} \log^2 n \right\}
\end{displaymath}
uniformly for $\pi \geq |\theta | \geq \delta_n$ and $n \geq n_0$.
\end{lemma}

\begin{proof} By taking logarithm, for $|x|<1$ and $|u|=1$, we get
\begin{displaymath}
\log G(u,x) = \log \left\{ \prod_{j=1}^\infty (1-ux^j)^{-j}
\right\} = - \sum_{j=1}^\infty j \log(1-u x^j)
\end{displaymath}
\begin{displaymath}
= \sum_{j=1}^\infty j\  \sum_{l=1}^\infty \frac{x^{jl} u^l}{l} =
\sum_{l=1}^\infty \frac{1}{l}\ \sum_{j=1}^\infty j (x^l)^j u^l.
\end{displaymath}
Thus, substituting $x=r_n e^{i\theta}$ and $u=e^{iT}$, we obtain
the estimate
\begin{displaymath}
|G(e^{iT}, r_n e^{i\theta})| = \left|\exp\left\{\log G(e^{iT}, r_n
e^{i\theta})\right\}\right| = \left|\ \exp\left\{\sum_{l=1}^\infty
\frac{1}{l}\ \sum_{j=1}^\infty j e^{ilT} r_n^{lj} e^{ijl\theta}
\right\}\right|
\end{displaymath}
\begin{displaymath}
=\left|\exp\left\{\sum_{j=1}^\infty j\ r_n^j\ e^{i(j\theta + T)} +
\sum_{l=2}^\infty \frac{1}{l} \ \sum_{j=1}^\infty j\ r_n^{lj}\
e^{il(j\theta + T)} \right\}\right|
\end{displaymath}
\begin{displaymath}
=\exp\left\{\sum_{j=1}^\infty j\ r_n^j\ \Re e\ e^{i(j\theta + T)}
+ \sum_{l=2}^\infty \frac{1}{l} \ \sum_{j=1}^\infty j\ r_n^{lj}\
\Re e\ e^{il(j\theta + T)} \right\}
\end{displaymath}
\begin{displaymath}
=\exp\left\{\sum_{j=1}^\infty j\ r_n^j\ \cos(j\theta + T) +
\sum_{l=2}^\infty \frac{1}{l}\ \sum_{j=1}^\infty j\ r_n^{lj}\
\cos[l(j\theta + T)] \right\}
\end{displaymath}
\begin{displaymath}
\leq\exp\left\{\sum_{j=1}^\infty j\ r_n^j\ \cos(j\theta + T) +
\sum_{l=2}^\infty\frac{1}{l}\ \sum_{j=1}^\infty j\ r_n^{lj}\
\right\}
\end{displaymath}
\begin{displaymath}
=\exp\left\{\sum_{j=1}^\infty j\ r_n^j\ [ \cos(j \theta + T)-1] +
\sum_{l=1}^\infty \frac{1}{l}\ \sum_{j=1}^\infty j\ r_n^{lj}\
\right\}
\end{displaymath}
\begin{equation}\label{eq.d1}
=G(1,r_n)\exp\left\{H_n(\theta, T)\right\},
\end{equation}
where
\begin{displaymath}
H_n(\theta, T)=\sum_{j=1}^\infty j\ r_n^j\ [\cos(j\theta + T)-1]
\end{displaymath}
\begin{equation}\label{eq.2.15}
=\Re e \left[\frac{r_n e^{i(\theta + T)}}{(1- r_n e^{i\theta})^2}
\right] - \frac{r_n}{(1-r_n)^2}
\end{equation}
\begin{displaymath}
= \frac{r_n\left[\cos(\theta + T) + r_n^2 \cos(\theta - T) - 2 r_n
\cos(T) \right]}{(1-2 r_n \cos(\theta) + r_n^2)^2} -
\frac{r_n}{(1-r_n)^2}.
\end{displaymath}
It is not difficult to show that the function
\begin{displaymath}
f(\theta)= \cos(\theta + T) + r_n^2 \cos(\theta - T)
\end{displaymath}
in the numerator of $H_n(\theta,T)$ attains its maximum value in
the range $\delta_n \leq \mid \theta \mid \leq \pi$ at $\theta =
\pm \delta_n$. To prove this one needs to determine the behaviour
of $f(\theta)$ in the ranges $T < \mit \theta \mit \leq \pi$ and
$\delta_n \leq \mit \theta \mit \leq T$. The first case is an easy
exercise whose study avoids the asymptotic form (\ref{def.T}) of
$T$. In the second case one can use the following representation
of the derivative
\begin{displaymath}
f^\prime (\theta)= -\sin(\theta + T) - r_n^2 \sin(\theta - T) =
-(1+ r_n^2) \theta -(1-r_n^2) T  + O(T^3).
\end{displaymath}
(\ref{eq.defrn}) and (\ref{def.T}) show that
\begin{displaymath}
(1 - r_n^2)T = \Theta(n^{-2/3}\log^{-1/2}{n}),\:\:  O(T^3) =
O(n^{-1} \log^{-3/2}{n}).
\end{displaymath}
Therefore, for sufficiently large $n$, the value of
$f^\prime(\theta)$ becomes close to $-(1 + r_n^2) \theta$. This
implies that $f(\theta)$ decreases for $\delta_n \leq \theta \leq
T$ and increases for $ -T  \leq \theta \leq -\delta_n$. It follows
that one may replace $\theta$ by $-\delta_n$ in (\ref{eq.2.15}).
Thus we can write
\begin{equation}\label{eq.hnt}
H_n(\theta, T) \leq \frac{r_n\left[\cos(-\delta_n + T) + r_n^2
\cos(-\delta_n - T) - 2 r_n \cos(T) \right]}{(1-2 r_n
\cos(\delta_n) + r_n^2)^2} - \frac{r_n}{(1-r_n)^2}.
\end{equation}
To estimate the cosine functions in the right-hand side of this
inequality we express $\delta_n$ and $T$ by (\ref{eq.defdelta})
and (\ref{def.T}), respectively. We get the following expansions:
\begin{eqnarray}
\cos \delta_n &=&1 - \frac{\delta_n^2}{2} + O(\delta_n^4) = 1 -
\frac{n^{-10/9}}{2\log^2{n}} + O(n^{-20/9} \log^{-4}{n}),\nonumber
\\
 \cos(-\delta_n \pm T)&=&\cos T \pm 2 \sin(T \mp \frac{\delta_n}{2})\, \sin \frac{\delta_n}{2},\nonumber  \\
 \cos T&=&1 - \Theta(n^{-2/3}/\log{n}).\nonumber
\end{eqnarray}
Moreover, (\ref{eq.defrn}) implies that
\begin{displaymath}
(1 - r_n)^{-2} = \left[\frac{n}{2 \zeta(3)} \right]^{2/3} +
O(n^{1/3}).
\end{displaymath}
Substituting these estimates in  the right-hand side of
(\ref{eq.hnt}), after some manipulations, we obtain
\begin{displaymath}
H_n(\theta, T) \leq \frac{r_n\left[ 1 - \Theta(n^{-2/3}/\log{n})
\right]}{(1-r_n)^2\left[1 + \frac{r_n n^{-10/9}}{(1-r_n)^2
\log^2{n}} +O(\frac{n^{-20/9}}{(1-r_n)^2 \log^4{n}}) \right]^2} -
\frac{r_n}{(1-r_n)^2}
\end{displaymath}
\begin{displaymath}
=\frac{r_n\left[ 1 - \Theta(n^{-2/3}/\log{n}) \right]}
{(1-r_n)^2\left[1 +
 2 n^{-4/9} / [2\zeta(3)]^{2/3} \log^2{n} +O(n^{-7/9} / \log^2{n}) \right]} -
\frac{r_n}{(1-r_n)^2}
\end{displaymath}
\begin{displaymath}
=\frac{r_n\left[ 1 - \Theta(n^{-2/3}/\log{n}) \right]} {(1-r_n)^2}
\left\{1 - \frac{2 n^{-4/9}} {[2\zeta(3)]^{2/3} \log^2{n}}
 +O\left(\frac{n^{-7/9}} {\log^2{n}}\right) \right\} - \frac{r_n}{(1-r_n)^2}
\end{displaymath}
\begin{displaymath}
 = -r_n \left\{ \left[ \frac{n}{2 \zeta(3)} \right]^{2/3} + O(n^{1/3}) \right\} \frac{2 n^{-4/9}} {[2\zeta(3)]^{2/3}
\log^2{n}} + O(1/\log{n})
\end{displaymath}
\begin{displaymath}
= -  \frac{2 n^{2/9}} {[2\zeta(3)]^{4/3} \log^2{n}} +
O(1/\log{n}). \nonumber
\end{displaymath}

Inserting this into (\ref{eq.d1}), we obtain the required bound.
\end{proof}


Further, we shall essentially use the asymptotic form of the
numbers $Q(n)$. It is given by Wright's formula (\ref{eq.wright}),
however, we need this result in a slightly different form. It is
not difficult to verify that Lemmas 1 and 2, (\ref{eq.defrn}) and
(\ref{eq.defb}) imply that $G(1, x)$ belongs to the class of
Hayman's admissible functions  \cite{Hayman} and thus, Hayman's
general asymptotic formula for the coefficients in the power
series representations of admissible functions is valid for the
number of plane partitions of $n$, $Q(n)$, as well. The next lemma
encompasses these results. We only sketch its proof and insert a
remark explaining the role of the asymptotic expansion
(\ref{eq.defrn}) there.

\begin{lemma}
For $G( 1,x )$, the generating function of the numbers $Q(n)$ of
plane partitions of $n$, defined by (\ref{brevg}), we have
\begin{equation}\label{eq.lemma3}
Q(n) \sim G( 1, r_n ) r_n^{-n} / [2 \pi b(r_n)]^{1/2}
\end{equation}
as $n\rightarrow \infty$, where $r_n$ satisfies the equation
\begin{equation}\label{eq.e1}
r G'( 1, r)/ G( 1, r) = n
\end{equation}
for sufficiently large $n$ and $b(r_n)$ is defined by
(\ref{eq.defb}).
\end{lemma}

\begin{sketch}
It is clear  that $|x|=1$ is a natural boundary for $G( 1, x)$.
Lemma 1 shows the behaviour of $G( 1, x)$ around its main
singularity $x=1$ (condition I of Hayman's Definition
\cite{Hayman}); Lemma 2 establishes the negligibility of the
growth of $G(1, x)$ as $x\rightarrow x_0$, $|x_0|=1$ and $x_0 \neq
1$ (condition II of \cite{Hayman}). It is then easily seen that
\begin{displaymath}
\frac{r\ G'(1,r)}{G(1,r)} = 2 \sum_{j=1} ^\infty
\frac{r^{2j}}{(1-r^j)^3}  + \sum_{j=1} ^\infty
\frac{r^j}{(1-r^j)^2}
\end{displaymath}
\begin{displaymath}
=\frac{2 \zeta(3)}{(1-r)^3} + o((1-r)^{-3})
\end{displaymath}
as $r\rightarrow 1^-$. This enables one to conclude that $r_n$,
determined for sufficiently large $n$ by (\ref{eq.e1}), can be
substituted by the asymptotic expansion (\ref{eq.defrn}). Thus one
can obtain (\ref{eq.lemma3}) after a direct application of
Hayman's theorem \cite{Hayman}.\hfill

Finally, we notice that (\ref{eq.defrn})-(\ref{eq.brn}) and
(\ref{eq.c1}) imply the coincidence of the right-hand sides of
(\ref{eq.wright}) and (\ref{eq.lemma3}).
\end{sketch}


\section{Proof of the Main Result}
\setcounter{equation}{0}

First, we let in (\ref{eq.stanley})
\begin{equation}\label{eq.deff}
G(u,x) = e^{F(u,x)},
\end{equation}
that is, we set
\begin{equation}\label{eq.f1}
F(u,x) = - \sum_{j=1}^\infty j \log(1-ux^j).
\end{equation}
We now apply Cauchy's coefficient formula to (\ref{eq.stanley}) on
the circle $x= r_n e^{i\theta}$, $-\pi <\theta \leq \pi$, with
$r_n$ determined by (\ref{eq.defrn}). Thus, in terms of the
notations (\ref{phi}), (\ref{eq.deff}) and (\ref{eq.f1}), for
$|u|\leq 1$, we obtain
\begin{displaymath}
Q(n) \varphi_n (u) = {r_n^{-n}\over 2\pi} \int_{-\pi}^\pi \exp
\left\{F(u,r_ne^{i\theta}) - i \theta n \right\} d\theta.
\end{displaymath}
We break up the range of integration as follows:
\begin{equation}\label{eq.int}
Q(n) \varphi_n(u) = J_1(n,u) + J_2(n,u),
\end{equation}
where
\begin{equation}\label{eq.defj1}
J_1(n,u) = {r_n^{-n}\over 2\pi} \int_{-\delta_n}^{\delta_n} \exp
\left\{F(u,r_ne^{i\theta}) - i \theta n \right\} d\theta,
\end{equation}
\begin{equation}\label{eq.defj2}
J_2(n,u) = {r_n^{-n}\over 2\pi} \int_{\delta_n\leq
|\theta|\leq\pi} \exp \left\{F(u,r_ne^{i\theta}) - i \theta n
\right\} d\theta.
\end{equation}

\subsection{An asymptotic estimate for $J_1(n,u)$.}

Using Taylor's formula expansion, we can write
\begin{displaymath}
F(u,r_n e^{i\theta})=F(u,r_n ) + r_n (e^{i\theta}-1) \left.
\frac{\partial}{\partial x} F(u,x) \right|
\begin{matrix} {\cr x=r_n}\end{matrix}
\end{displaymath}
\begin{equation}\label{eq.g1}
+ \frac{r_n^2}{2} (e^{i\theta}-1)^2
\left.\frac{\partial^2}{\partial x^2}F(u,x) \right|
\begin{matrix} {\cr x=r_n}\end{matrix} + O\left(|\theta|^3  \left.
\frac{\partial^3}{\partial x^3}F(|u|,x)\right|\begin{matrix}{\cr
x=r_n}\end{matrix} \right).
\end{equation}
To evaluate the partial derivatives of $F(u,x)$, we shall use the
following expressions:
\begin{equation}\label{eq.g2}
\left.\frac{\partial}{\partial x}F(u,x)\right|\begin{matrix}{\cr
x=r_n}\end{matrix} = u \sum_{j=1}^{\infty} \frac{j^2
r_n^{j-1}}{1-u r_n^j},
\end{equation}
\begin{equation}\label{eq.g3}
\left.\frac{\partial^2} {\partial x^2} F(u,x ) \right|
\begin{matrix} {\cr x=r_n} \end{matrix} = u \sum_{j=2}^{\infty}
\frac{j^2 (j-1) r_n^{j-2}} {1-u r_n^j} + u^2 \sum_{j=1}^{\infty}
\frac{j^3 r_n^{2(j-1)}}{(1-u r_n^j)^2},
\end{equation}
\begin{displaymath}
\left.\frac{\partial^3} {\partial x^3} F(u,x ) \right|
\begin{matrix} {\cr x=r_n}\end{matrix} = u \sum_{j=3}^{\infty}
\frac{j^3 (j-1) (j-2)r_n^{j-3}} {1-u r_n^j}
\end{displaymath}
\begin{equation}\label{eq.g4}
+ 3u^2 \sum_{j=2}^{\infty} \frac{j^3(j-1)r_n^{2j-3}}{(1-u
r_n^j)^2}+ 2u^3 \sum_{j=1}^{\infty} \frac{j^4 r_n^{3(j-1)}} {(1-u
r_n^j)^3}.
\end{equation}

We proceed further to establish the convergence of $J_1(n,u)$ in
terms of Fourier transforms by setting $u=e^{iT}$,
$-\infty<T<\infty$, in (\ref{eq.defj1}),
(\ref{eq.g1})-(\ref{eq.g4}). In what follows later, an application
of L\'{e}vy's continuity theorem for characteristic functions
\cite[Section 3.6]{Lukacz} would specify the value of $T$ as a
function of the main parameter $n$. Furthermore, we need notations
for the following four functions:
\begin{equation}\label{eq.defpsi}
\psi_{m,k}(z,T) = \int_z^\infty \frac{y^m}{(e^{y-iT}-1)^{m-k}}dy,
\end{equation}
\begin{displaymath}
(k,m)= (0,1), (1,2), (1,3), (2,3),\; 0\leq z < \infty.
\end{displaymath}
They are closely related to the Debye functions (see e.g.
\cite[Section 27.1]{Abramovitz}). Moreover, formula 27.1.3 of
\cite{Abramovitz} shows that
\begin{equation}\label{eq.g5}
\psi_{m,m-1}(0,0) = m! \zeta(m+1).
\end{equation}

We now proceed to the asymptotic estimates of the summands in
(\ref{eq.g1}). Interpreting again the sum in (\ref{eq.g2}) by a
Riemann's one with the same step size $y_n = - \log r_n$ as in the
proof of Lemma 2 of \cite{Mutafchiev2} and using (\ref{eq.g5}) and
(\ref{eq.defyn}), we find that
\begin{displaymath}
r_n \,\left.\frac{\partial}{\partial
x}F(u,x)\right|\begin{matrix}{\cr x=r_n, u=e^{iT} }\end{matrix} =
\psi_{2,1} (y_n,T)y_n^{-3} + O(y_n^{-1})
\end{displaymath}
\begin{displaymath}
=[\psi_{2,1} (0,T) + O(n^{-2/3})]y_n^{-3} + O(n^{1/3})
\end{displaymath}
\begin{displaymath}
=[\psi_{2,1} (0,0) + R(n,T)]y_n^{-3} + O(n^{1/3})
\end{displaymath}
\begin{equation}\label{firstdF}
=[2\zeta(3) + R(n,T)]\left[\frac{n}{2\zeta(3)}+ O(1)\right] +
O(n^{1/3})
\end{equation}
\begin{displaymath}
= n + \frac{n\, R(n,T)}{2\,\zeta(3)} + O(| R(n,T)|) + O( n^{1/3}),
\end{displaymath}
where
\begin{displaymath}
R(n,T) = T\, \psi_{2,1}^\prime (0, T_1) = i\, T \, \int_{0}^\infty
\frac{y^2\,e^{y- i\, T_1}\, dy}{(e^{y- i\, T_1}-1)^2}
\end{displaymath}
is the remainder term in the Taylor's formula expansion of
$\psi_{2,1}(0, T)$ and $0< T_1 <T$. Since the last integral is
finite, from (\ref{def.T}) we get
\begin{equation}\label{rnt}
R(n,T) = O(n^{-1/3} / \sqrt{\log n}).
\end{equation}
Furthermore, note that (\ref{eq.g1}) requires a multiplication by
$e^{i \theta }-1$. Hence by (\ref{eq.defrn}), (\ref{eq.defdelta})
(\ref{firstdF}) (\ref{rnt})
\begin{equation}\label{eq.g6}
r_n(e^{i\theta}-1)\left.\frac{\partial}{\partial
x}F(u,x)\right|\begin{matrix}{\cr x=r_n, u=e^{iT}} \end{matrix}
\end{equation}
\begin{displaymath}
= [i\theta + O(|\theta|^2)]\left[n + \frac{n\,
R(n,T)}{2\,\zeta(3)} + O( n^{1/3})\right]
\end{displaymath}
\begin{displaymath}
= i\theta n + i\, \theta \,\frac{n\, R(n,T)}{2\,\zeta(3)} +
O(\delta_n\, n^{1/3}) +O(n\delta_n^2)
\end{displaymath}
\begin{displaymath}
=i\theta n +i\, \theta \,\frac{n\, R(n,T)}{2\,\zeta(3)} +
O(n^{-1/9}/\log^2 n).
\end{displaymath}
To deal with the third term of the expansion in (\ref{eq.g1}), one
has to follow the same line of reasoning. Thus, in a similar way
(\ref{eq.g3}) becomes
\begin{equation}\label{eq.g7}
r_n^2\,\left.\frac{\partial^2}{\partial x^2}F(u,x)\right|_{ x=r_n,
u=e^{iT} } = y_n^{-4}[\psi_{3,2}(y_n,T) + \psi_{3,1}(y_n,T)]
+O(n).
\end{equation}
On the other side, using Taylor's formula expansion and
(\ref{eq.g5}), we obtain
\begin{equation}\label{eq.g8}
\psi_{3,2}(y_n,T) = \psi_{3,2}(0,T) +
O\left(y_n\left|\left.\frac{\partial}{\partial z} \psi_{3,2}(z,T)
\right|_{z=y_n}\right|\right)
\end{equation}
\begin{displaymath}
= \psi_{3,2}(0,T) + O(1/n) =  \psi_{3,2}(0,0) + O(|T|) +O(1/n)=
6\zeta(4) + O(n^{-1/3} / \sqrt{\log n}),
\end{displaymath}
\begin{equation}\label{eq.g9}
\psi_{3,1}(y_n,T) = \psi_{3,1}(0,T) +
O\left(y_n\left|\left.\frac{\partial}{\partial z} \psi_{3,1}(z,T)
\right|_{z=y_n}\right|\right)
\end{equation}
\begin{displaymath}
= \psi_{3,1}(0,T) + O(n^{-2/3}) =  \psi_{3,1}(0,0) + O(|T|)
+O(n^{-2/3})
\end{displaymath}
\begin{displaymath}
= 6[\zeta(3) -\zeta(4)] + O(n^{-1/3} / \sqrt{\log n}).
\end{displaymath}
So (\ref{eq.defrn}), (\ref{eq.defdelta}), (\ref{eq.defyn}) and
(\ref{eq.g7}) - (\ref{eq.g9}) imply
\begin{equation}\label{eq.g10}
\frac{r_n^2}{2}\,(e^{i\theta}-1)^2\left.\frac{\partial^2}{\partial
x^2}F(u,x)\right|_{x=r_n, u=e^{iT} }
\end{equation}
\begin{displaymath}
=\frac{1}{2}[1+O(n^{-1/3})][-\theta^2 +
O(|\theta|^3)]\left[\left(\frac{n}{2\zeta(3)} \right)^{4/3}
+O(n^{1/3})\right]
\end{displaymath}
\begin{displaymath}
\times \{ 6 \zeta(4)+ 6 [\zeta(3) - \zeta(4)] +  O(n^{-1/3} /
\sqrt{\log n}) \}
\end{displaymath}
\begin{displaymath}
=\frac{1}{2}[1+O(n^{-1/3})]  [-\theta^2 + O(\delta_n^3)] \left[
\left( \frac{n}{2\zeta(3)} \right)^{4/3}+O(n^{1/3})\right]
\end{displaymath}
\begin{displaymath}
 \times [6\zeta(3)+  O(n^{-1/3} / \sqrt{\log n})]
\end{displaymath}
\begin{displaymath}
=\frac{-\theta^2}{2} \left( \frac{n}{2\zeta(3)} \right)^{4/3}
[6\zeta(3)][1 + O(n^{-1/3})] + O(n^{1/3} \delta_n^2) +O(n^{4/3}
\delta_n^3 )
\end{displaymath}
\begin{displaymath}
=\frac{-\theta^2}{2} \left( \frac{n}{2\zeta(3)} \right)^{4/3}
[6\zeta(3)] + O(n^{4/3} \delta_n^3) +  O(n\delta_n^2)
\end{displaymath}
\begin{displaymath}
=\frac{-\theta^2}{2} \left( \frac{n}{2\zeta(3)} \right)^{4/3}
[6\zeta(3)] +  O(n^{-1/9}/\log^2 n).
\end{displaymath}
Finally, following similar but simpler analysis as above, with the
aid of (\ref{eq.g2}) we can show the negligibility of the error
term in (\ref{eq.g1}). We have
\begin{displaymath}
O\left(|\theta|^3  \left. \frac{\partial^3}{\partial
x^3}F(|u|,x)\right|_{ x=r_n} \right) = O(\delta_n^3 y_n^{-5})
\end{displaymath}
\begin{equation}\label{eq.g11}
= O((n^{-15/9}/\log^3 n) n^{5/3}) = O(1/\log^3 n).
\end{equation}

We are now ready to substitute (\ref{eq.g1}), (\ref{eq.g6}),
(\ref{eq.g10}) and (\ref{eq.g11}) into the integral of
$J_1(n,e^{iT})$ (see the (\ref{eq.defj1})). We can write
\begin{displaymath}
J_1(n,e^{iT}) = \frac{r_n^{-n}e^{F(e^{iT}, r_n)}}{2\pi}
\int_{-\delta_n}^{\delta_n} \exp \left\{ i\theta n +
i\theta\frac{n R(n,T)}{2\,\zeta(3)} \right.
\end{displaymath}
\begin{displaymath}
\left. -\frac{\theta^2}{2}\left[\frac{n}{2\zeta(3)}\right]^{4/3}
[6\zeta(3)] +  O(n^{-1/9}/\log^2 n) + O(1/\log^3 n) -i\theta n
\right\} d\theta
\end{displaymath}
\begin{equation}\label{eq.g12}
=\frac{r_n^{-n}e^{F(e^{iT}, r_n)}}{2\pi b^{1/2}(r_n)} \left[ 1 +
O\left(\frac{1}{\log^3 n}\right) \right]
\int_{-\delta_n  b^{1/2}(r_n)}^{\delta_n b^{1/2}(r_n)} \exp
\left\{ \frac{i t n^{1/3} \,R(n,T)}{[2\,\zeta(3)]^{5/6}\sqrt{3}} -
\frac{t^2}{2}\right\}dt.
\end{equation}
Note that we substituted $\theta = t/  b^{1/2}(r_n), -\infty < t <
\infty$, with $ b^{1/2}(r_n)$ defined by (\ref{eq.defb}) and
(\ref{eq.brn}). In addition, (\ref{eq.brn}) and
(\ref{eq.defdelta}) justify the computation
\begin{equation}\label{eq.g13}
\delta_n  b^{1/2}(r_n) \sim d n^{1/9} \log n,\; d= \sqrt{3} /
[2\zeta(3)]^{1/6}.
\end{equation}
So, the bounds of the last integral in (\ref{eq.g12}) tend  to
$\pm \infty$. The estimate (\ref{rnt}) for $R(n,T)$ and the
Lebesgue's dominated convergence theorem allow the passage to the
limit under the integral. Moreover, (\ref{rnt}) shows that the
limit of the integrand equals $e^{-t^2/2}$ with an error term of
order $O(1/\sqrt{\log n})$. The additional error term that we get
replacing $\pm \delta_n  b^{1/2}(r_n)$ by $\pm \infty$ can be
estimated using the asymptotic expansion of the function
$(2\pi)^{-1/2} \int_{z}^{\infty} e^{-t^2} dt $ as $z\rightarrow
\infty$ (see \cite[Chap.7] {Abramovitz}). Therefore,
(\ref{eq.g13}) implies that this error term is at most $O(n^{-1/9}
(\log n) \exp\{-d^2 n^{2/9} \log^2 n\})= O(1/\sqrt{\log n})$. This
shows in turn that
\begin{displaymath}
J_1(n,e^{iT}) =\frac{r_n^{-n}e^{F(e^{iT}, r_n)}}{2\pi
b^{1/2}(r_n)} \left[ 1 + O(1/\sqrt{\log n}) \right]
\int_{-\infty}^{\infty} e^{ - \frac{t^2}{2}} \, dt
\end{displaymath}
\begin{equation}\label{eq.g14}
=\frac{r_n^{-n}e^{F(e^{iT}, r_n)}}{\sqrt{2\pi b^(r_n)}}\left[ 1 +
O(1/\sqrt{\log n}) \right].
\end{equation}
This completes the required estimate for $J_1(n, e^{iT})$.

\subsection{An asymptotic estimate for $J_2(n,u)$.}

To estimate this integral asymptotically we only need to apply
directly Lemmas 2 and 3. So for $n\geq n_0$ and bounded $T$, we
get from (\ref{eq.defj1}) and (\ref{eq.brn}):
\begin{displaymath}
|J_2(n,u)| \leq \frac{r_n^{-n}}{2\pi} \int_{\pi \geq |\theta| \geq
\delta_n} \left|G(e^{iT}, r_n e^{i\theta})\right|d\theta
\end{displaymath}
\begin{equation}\label{eq.h1}
\leq r_n^{-n} G(1, r_n)(1+ \varepsilon) \exp\{-2n^{2/9}/[2
\zeta(3)]^{4/3} \log^2 n\}
\end{equation}
\begin{displaymath}
=\frac{r_n^{-n} G(1, r_n)}{\sqrt{2\pi b(r_n)}} (1+ \varepsilon)
\exp\{-2n^{2/9}/[2 \zeta(3)]^{4/3} \log^2 n\ + \log \sqrt{2\pi
b(r_n)}\}
\end{displaymath}
\begin{displaymath}
\sim Q(n)(1+ \varepsilon) \exp\{-2n^{2/9}/[2 \zeta(3)]^{4/3}
\log^2 n\ +\frac{2}{3} \log n +O(1)\} = o(Q(n)).
\end{displaymath}

\subsection{Final estimate for $J_1(n,e^{iT})$: the choice of $T=T(n)$}

We start our estimation with a Taylor's formula expansion:
\begin{displaymath}
F(e^{iT}, r_n) = F(1, r_n) + (e^{iT}-1) \left.
\frac{\partial}{\partial u} F(u, r_n) \right|_{u=1}
\end{displaymath}
\begin{displaymath}
+ (e^{iT}-1)^2 \left. \frac{\partial^2}{\partial u^2} F(u,
r_n)\right|_{u=1} +O\left(|T|^3 \left.\left|
\frac{\partial^3}{\partial u^3} F(u, r_n) \right|_{u=1} \right|
\right)
\end{displaymath}
\begin{equation}\label{eq.i1}
= F(1, r_n) + iT \left. \frac{\partial}{\partial u} F(u,
r_n)\right|_{u=1} - T^2 \left. \frac{\partial^2}{\partial u^2}
F(u, r_n)\right|_{u=1}
\end{equation}
\begin{displaymath}
+ O\left(|T|^3 \left[ \left| \left. \frac{\partial^3}{\partial
u^3} F(u, r_n) \right|_{u=1} \right| + \left| \left.
\frac{\partial^2}{\partial u^2} F(u, r_n) \right|_{u=1} \right|
\right] + |T|^2 \left| \left. \frac{\partial}{\partial u} F(u,
r_n) \right|_{u=1} \right| \right).
\end{displaymath}
The partial derivatives of $F(u, r_n)$ can be again evaluated by
Riemann's integral in a way similar to that presented in the proof
of Lemma 2 of \cite{Mutafchiev2}. We  have by (\ref{eq.defyn}),
(\ref{eq.defpsi}) and (\ref{eq.g5}) that
\begin{displaymath}
\left. \frac{\partial}{\partial u} F(u, r_n)\right|_{u=1} =
\sum_{j=1}^\infty \frac{jr_n^j}{1-r_n^j} = y_n^{-2}
\int_{y_n}^\infty \frac{v e^{-v}}{1- e^{-v}} dv + O(y_n^{-1})
\end{displaymath}
\begin{displaymath}
= \left[ \frac{n}{2 \zeta(3)} \right]^{2/3} \left[ 1 + O(n^{-1})
\right] \psi_{1, 0}(0,0) \left[ 1 + O(n^{-1/3}) \right]
+O(n^{1/3})
\end{displaymath}
\begin{equation}\label{eq.i2}
=\left[ \frac{n}{2 \zeta(3)} \right]^{2/3} \zeta(2) + O(n^{1/3}).
\end{equation}
Similarly
\begin{equation}\label{eq.i3}
\left. \frac{\partial^2}{\partial u^2} F(u, r_n)\right|_{u=1} =
\sum_{j=1}^\infty \frac{jr_n^{2j}} {(1-r_n^j)^2} = y_n^{-2}
I(y_n)+ o(y_n^2 I(y_n)),
\end{equation}
where
\begin{displaymath}
I(y_n) =  \int_{y_n}^\infty \frac{v e^{-2v}}{(1- e^{-v})^2} dv =
\frac{y_n}{e^{y_n}-1} - \psi_{1,0}(y_n,0) + \int_{y_n}^\infty
\frac{dv}{e^{v}-1}.
\end{displaymath}
Since the first two summands above have finite limits and
\begin{displaymath}
\frac{ \int_{y_n}^\infty \frac{dv}{e^{v}-1}}{-\log y_n}= 1+ o(1)
\end{displaymath}
as $y_n\rightarrow 0$, we observe that
\begin{displaymath}
I(y_n) = O(1) - \log y_n [1+ o(1)] =  [1+ o(1)](- \log y_n).
\end{displaymath}
Substituting this into (\ref{eq.i3}) and invoking the asymptotic
(\ref{eq.defyn}) for $y_n$, we find that
\begin{equation}\label{eq.i4}
\left. \frac{\partial^2}{\partial u^2} F(u, r_n)\right|_{u=1} =
\left[ \frac{n}{2 \zeta(3)}\right]^{2/3} \frac{\log n}{3} + o
\left(n^{2/3}\log n\right).
\end{equation}
In the same way one can deduce the following estimate for the
third partial derivative:
\begin{equation}\label{eq.i5}
\left. \frac{\partial^3}{\partial u^3} F(u, r_n)\right|_{u=1} =
O(n).
\end{equation}
Putting (\ref{eq.i2}), (\ref{eq.i4}) and (\ref{eq.i5}) into
(\ref{eq.i1}), we obtain
\begin{displaymath}
F(e^{iT}, r_n) = F(1, r_n) + i T  \left[ \frac{n}{2
\zeta(3)}\right]^{2/3} \zeta(2) +O(|T|) -\frac{T^2}{2} \left[
\frac{n}{2 \zeta(3)}\right]^{2/3} \frac{\log n}{3}
\end{displaymath}
\begin{equation}\label{eq.i6}
+ o \left(T^2 n^{2/3}\log n\right) + O(n|T|^3) +O(n^{2/3}T^2).
\end{equation}

It is now clear how to choose the variable $T=T(n)$ which
satisfies (\ref{def.T}). Note that this choice determines the
scaling factor in our limit theorem. Setting
\begin{equation}\label{eq.i7}
T= w\left/  \left[ \frac{n}{2 \zeta(3)}\right]^{1/3}\right.
\sqrt{\frac{ \log n}{3}}, \; -\infty < w < \infty,
\end{equation}
(in agreement with (\ref{def.T})) we can rewrite (\ref{eq.i6}) in
the following form
\begin{equation}\label{eq.i8}
F\left(\exp \left\{iw\left/  \left[ \frac{n}{2
\zeta(3)}\right]^{1/3}\right. \sqrt{\frac{\log n}{3}} \right\},
r_n\right)
\end{equation}
\begin{displaymath}
= F(1, r_n) + i w  \left[ \frac{n}{2 \zeta(3)}\right]^{1/3}
\zeta(2)\left/ \sqrt{\frac{\log n}{3}} - \frac{w^2}{2} +
o(1)\right. .
\end{displaymath}
The convergence here is uniform with respect to $w$ belonging to
any finite interval. Going back to the estimate (\ref{eq.g14}) for
$J_1(n,e^{iT})$ and replacing there (\ref{eq.i7}) and
(\ref{eq.i8}), we find that
\begin{equation}\label{eq.i9}
J_1\left(n, \exp \left\{iw\left/  \left[ \frac{n}{2
\zeta(3)}\right]^{1/3}\right. \sqrt{\frac{\log n}{3}} \right\}
\right)=\frac{r_n^{-n}e^{F(1, r_n)}}{\sqrt{2\pi b^(r_n)}}
\end{equation}
\begin{displaymath}
\times \left[1+ O(\frac{1}{\sqrt{\log n}})\right]  \exp \left\{ i
w \left[ \frac{n}{2 \zeta(3)}\right]^{1/3} \zeta(2) \left/
\sqrt{\frac{\log n}{3}}\right. - \frac{w^2}{2} + o(1) \right\}
\end{displaymath}
\begin{displaymath}
\sim Q(n)\exp \left\{ i w \left[ \frac{n}{2 \zeta(3)}\right]^{1/3}
\zeta(2) \left/ \sqrt{ \frac{\log n}{3}}\right. - \frac{w^2}{2}
\right\},
\end{displaymath}
where the reduction in the last asymptotic equivalence is due to a
direct application of Lemma 3. To handle with $J_1$ in a more
appropriate way we may rewrite (\ref{eq.i9}) as follows:
\begin{displaymath}
J_1\left(n, \exp \left\{i w \left/  \left[ \frac{n}{2
\zeta(3)}\right]^{1/3}\right. \sqrt{\frac{\log n}{3}} \right\}
\right)
\end{displaymath}
\begin{equation}\label{eq.i10}
\times \exp \left\{- i w \left[ \frac{n}{2 \zeta(3)}\right]^{1/3}
\zeta(2)\left/ \sqrt{\frac{\log n}{3}} \right. \right\} \sim Q(n)
e^{-w^2/2}.
\end{equation}
It remains now to deal with the characteristic function of
$(\tau_n - c_0 n^{2/3})/ c_1 n^{1/3} \log^{1/2} n$. The
expressions for $c_0$ and $c_1$ in the Theorem and
(\ref{eq.stanley}) imply that it is equal to
\begin{displaymath}
\exp \left\{- i w \left[ \frac{n}{2 \zeta(3)}\right]^{1/3}
\zeta(2) \left/ \sqrt{\frac{\log n}{3}}\right. \right\}
\end{displaymath}
\begin{displaymath}
\times \varphi_n \left(\exp \left\{i w \left/  \left[ \frac{n}{2
\zeta(3)}\right]^{1/3}\right. \sqrt{\frac{\log n}{3}} \right\}
\right)
\end{displaymath}
\begin{displaymath}
= \exp\left\{ -i w \frac{c_0 n^{1/3}}{c_1 \log^{1/2}n} \right\}
\varphi_n(\exp\{ i w / c_1 n^{1/3} \log^{1/2}n\}) = \Phi_n(w)
\end{displaymath}
Setting $u=e^{iT}$ in (\ref{eq.int}) with $T$ given by
(\ref{eq.i7}) and substituting estimates (\ref{eq.i9}) and
(\ref{eq.h1}) in it, we obtain
\begin{displaymath}
Q(n) \Phi_n(w) = Q(n) e^{-w^2/2} + o(Q(n)).
\end{displaymath}
The required week convergence follows from L\`{e}vy's continuity
theorem for characteristic functions \cite[Section 3.6]{Lukacz}.
\vskip0,5cm

{\bf Acknowledgements} \vspace{,5cm}

We thank Cyril Banderier who kindly informed us about the results
on plane partitions contained in \cite{AlmkvistII}. We are also
grateful to the referee for the careful reading of the manuscript
and especially for his help to eliminate defects in the choice
(\ref{eq.defrn}) of $r_n$ and in the estimate of the sum in
(\ref{eq.g2}).


\vspace{2cm}

Ljuben R. Mutafchiev

American University in Bulgaria

2700 Blagoevgrad

Bulgaria

email: ljuben@aubg.bg

\vspace{.5cm}

Emil P. Kamenov

Sofia University "St. Kliment Ohridski"

Faculty of Mathematics and Informatics

James Bouchier Blvd. 5

Sofia 1164

Bulgaria

email: kamenov@fmi.uni-sofia.bg

\end{document}